\documentclass[12pt]{amsart}
\usepackage{amssymb,amsmath,latexsym}

\newcommand\eset{\varnothing}
\newcommand\bbP{\mathbb{P}} 
\newcommand\cA{{\mathcal A}}
\newcommand\cC{{\mathcal C}}
\newcommand\cL{{\mathcal L}}
\newcommand\cM{{\mathcal M}}
\newcommand\cP{{\mathcal P}}
\newcommand\gauss[2]{\left[ { #1 \atop #2 } \right]} 
\newcommand\smallgauss[2]{\big[ { #1 \atop #2 } \big]}
\renewcommand\qedsymbol{$\blacksquare$} 

\newtheorem{theorem}{Theorem} 
\newtheorem{lemma}[theorem]{Lemma} 
\newtheorem{corollary}[theorem]{Corollary} 
\newtheorem{proposition}[theorem]{Proposition} 

\begin{document}

\begin{center}

\Large{\bf A Meshalkin Theorem for Projective Geometries}\footnote{Appeared in \emph{Journal of Combinatorial Theory Series A} {\bf 102} (2003), 433-441.}  

\normalsize

{\sc Matthias Beck and Thomas Zaslavsky\footnote{Research supported by National Science Foundation grant DMS-0070729.}}

\bigskip

{\sc Department of Mathematical Sciences \\
     State University of New York at Binghamton \\
     Binghamton, NY, U.S.A.\ 13902-6000}

\bigskip

{\tt matthias@math.binghamton.edu \\
     zaslav@math.binghamton.edu }

\bigskip\bigskip

Dedicated to the memory of Lev Meshalkin.

\end{center}

\bigskip\bigskip\bigskip\bigskip 

%\footnotesize
{\it Abstract:} 
%\begin{abstract}
Let $\cM$ be a family of sequences $ ( a_1 , \dots, a_p ) $ where each $a_k$ is a flat in a 
projective geometry of rank $n$ (dimension $n-1$) and order $q$, and the sum of ranks, 
$ r(a_1) + \dots + r(a_p) $, equals the rank of the join $ a_1 \vee \dots \vee a_p $. 
We prove upper bounds on $|\cM|$ and corresponding LYM inequalities assuming that 
(i) all joins are the whole geometry and for each $k<p$ the set of all $a_k$'s of 
sequences in $\cM$ contains no chain of length $l$, and that 
(ii) the joins are arbitrary and the chain condition holds for all $k$. 
These results are $q$-analogs of generalizations of Meshalkin's and Erd\H{o}s's 
generalizations of Sperner's theorem and their LYM companions, and they generalize Rota 
and Harper's $q$-analog of Erd\H{o}s's generalization. 
%\end{abstract}
\normalsize



%{Extremal set theory, Combinatorial structures in finite projective spaces, Combinatorics of posets.}

%\emph{Running head}: Projective Meshalkin theorem

\vfill\pagebreak
% -------------------------------------------

\section{Introducing the Players}\label{intro} 

We present a theorem that is at once a $q$-analog of a generalization, due to Meshalkin, of Sperner's famous theorem on antichains of sets and a generalization of Rota and Harper's $q$-analog of both Sperner's theorem and Erd\H{o}s's generalization. 

Sperner's theorem \cite{sperner} concerns a subset $\cA$ of $\cP(S)$, the power set of an $n$-element set $S$, that is an \emph{antichain}: no member of $\cA$ contains another.  
It is part (b) of the following theorem. Part (a), which easily implies (b) (see, e.g., \cite[Section 1.2]{anderson}) was found later by Lubell \cite{lubell}, Yamamoto \cite{yamamoto}, and Meshalkin \cite{meshalkin} 
(and Bollob\'as independently proved a generalization \cite{bollobas}); consequently, it and similar inequalities are called 
\emph{LYM inequalities}. 

\begin{theorem}\label{s} 
Let $\cA$ be an antichain of subsets of $S$. Then: 
\begin{enumerate} 
\item[(a)] $\displaystyle \sum_{A\in \cA} \frac {1}{|A|} \leq 1$ and 
\item[(b)] $|\cA | \leq {n\choose {\lfloor n/2 \rfloor }}$. 
\item[(c)] Equality occurs in (a) and (b) if $\cA$ consists of all subsets of $S$ of size $ \lfloor n/2 \rfloor $, or all of size $\lceil n/2 \rceil$. 
\end{enumerate} 
\end{theorem} 

The idea of Meshalkin's insufficiently well known generalization\footnote{We do not 
find it in books on the subject \cite{anderson,engel} but only in \cite{klain}.} 
(an idea he attributes to Sevast'yanov) is to consider ordered $p$-tuples 
$A=(A_1,\dots,A_p)$ of pairwise disjoint sets whose union is $S$. We call these 
\emph{weak compositions of $S$ into $p$ parts}. 

\begin{theorem}\label{m} 
Let $\cM$ be a family of weak compositions of $S$ into $p$ 
parts such that each set $\cM_k = \{ A_k : \ A \in \cM \}$ is an antichain. 
\begin{enumerate} 
\item[(a)] $ \displaystyle \sum_{A\in \cM} \dfrac{1}{\binom{n}{|A_1|,\hdots,|A_p|}} \leq 1 $. 
\item[(b)] $ \displaystyle |\cM| \leq \max_{ \alpha_1 + \dots + \alpha_p = n } \binom{n}{\alpha_1,\dots,\alpha_p} = \binom n { \big\lceil \frac n p \big\rceil , \dots , \big\lceil \frac n p \big\rceil , \big\lfloor \frac n p \big\rfloor , \dots , \big\lfloor \frac n p \big\rfloor } $. 
\item[(c)] Equality occurs in (a) and (b) if, for each $k$, $\cM_k$ consists of all subsets of $S$ of size $ \big\lceil \frac n p \big\rceil $, or all of size $\big\lfloor \frac n p \big\rfloor$. 
\end{enumerate} 
\end{theorem} 

Part (b) is Meshalkin's theorem \cite{meshalkin}; the corresponding LYM inequality (a) was subsequently found by Hochberg and Hirsch \cite{hh}. 
(In expressions like the multinomial coefficient in (b), since the lower numbers must sum to $n$, 
the number of them that equal $ \big\lceil \frac n {p} \big\rceil $ is the least nonnegative residue of $n$ modulo $p+1$.)

In \cite{bwz} Wang and we generalized Theorem \ref{m} in a way that simultaneously also generalizes Erd\H{o}s's theorem on \emph{$l$-chain-free families}: subsets of $\cP(S)$ that contain no chain of length $l$. (Such families have been called ``$r$-families'' and ``$k$-families'', where $r$ or $k$ is the forbidden length. We believe a more suggestive name is needed.) 

\begin{theorem}[{\cite[Corollary 4.1]{bwz}}]\label{bwz} 
Let $\cM$ be a family of weak compositions of $S$ into $p$ parts such that each $\cM_k$, for $k<p$, is $l$-chain-free.  Then:
\begin{enumerate} 
\item[(a)] $ \displaystyle \sum_{A\in \cM} \dfrac{1}{\binom{n}{|A_1|,\hdots,|A_p|}} \leq l^{p-1} $, and 
\item[(b)] $|\cM|$ is no greater than the sum of the $l^{p-1}$ largest multinomial coefficients of the form $ \binom{n}{\alpha_1,\dots,\alpha_p} $. 
\end{enumerate} 
\end{theorem} 

Erd\H{o}s's theorem \cite{erdos} is essentially the case $p=2$, in which $A_2 = S \setminus A_1$ is redundant. 
The upper bound is then the sum of the $l$ largest binomial coefficients $ \binom n j , \ 0 \leq j \leq n $, 
and is attained by taking a suitable subclass of $\cP(S)$. In general the bounds in Theorem 
\ref{bwz} cannot be attained \cite[Section 5]{bwz}. 

Rota and Harper began the process of $q$-analogizing by finding versions 
of Sperner's and Erd\H{o}s's theorems for finite projective geometries \cite{rota}. 
We think of a projective geometry $ \bbP^{\,n-1} = \bbP^{\,n-1} (q) $ of order $q$ and 
rank $n$ (i.e., dimension $n-1$) as a lattice of flats, in which $ \hat 0 = \eset $ 
and $\hat 1$ is the whole set of points. The \emph{rank} of a flat $a$ is $r(a)= \dim a + 1 $. 
The \emph{$q$-Gaussian coefficients} (usually the ``$q$'' is omitted) are the quantities 
\[ 
\gauss n k = \frac{ n !_q }{ k!_q (n-k) !_q } \qquad \text{ where } \qquad n !_q = ( q^n - 1 ) ( q^{n-1} - 1 ) \cdots (q-1) \ . 
\] 
They are the $q$-analogs of the binomial coefficients. Again, a family of projective flats 
is \emph{$l$-chain-free} if it contains no chain of length $l$.  
Let $\cL_k$ be the set of all flats of rank $k$ in $\bbP^{n-1}(q)$.

\begin{theorem}[{\cite[p.~200]{rota}}]\label{rh} 
Let $\cA$ be an $l$-chain-free family 
of flats in $ \bbP^{\,n-1} (q) $. 
\begin{enumerate} 
\item[(a)] $ \displaystyle \sum_{a \in \cA} \dfrac{1}{\smallgauss{n}{r(a)}} \leq l $. 
\item[(b)] $|\cA|$ is at most the sum of the $l$ largest Gaussian coefficients $ \smallgauss n j $ for $ 0 \leq j \leq n $. 
\item[(c)] There is equality in (a) and (b) when $\cA$ consists of the $l$ largest classes $\cL_k$, if $n-l$ is even, or the $l-1$ largest 
classes and one of the two next largest classes, if $n-l$ is odd. 
\end{enumerate} 
\end{theorem} 

Our $q$-analog theorem concerns the projective analogs of weak compositions of a set. 
A \emph{Meshalkin sequence of length $p$} in $ \bbP^{\,n-1} (q) $ is a sequence 
$a = (a_1,\dots,a_p)$ of flats whose join is $\hat 1$ and whose ranks sum to $n$. 
The submodular law implies that, if $a_J := \bigvee_{j \in J} a_j$ 
for an index subset $J \subseteq [p] = \{ 1,2,\dots,p \} $, then $ a_I \wedge a_J = \hat 0 $ 
for any disjoint $ I, J \subseteq [p] $; so the members of a Meshalkin sequence are highly disjoint. 

To state the result we need a few more definitions. If $\cM$ is a set of 
Meshalkin sequences, then for each $ k \in [p] $ we define 
$ \cM_k := \{ a_k : \ (a_1,\dots,a_p) \in \cM \} $. If $ \alpha_1,\dots,\alpha_p $ 
are nonnegative integers whose sum is $n$, we define the \emph{($q$-)Gaussian multinomial coefficient} 
to be 
\[ 
\gauss n \alpha = \gauss n { \alpha_1,\dots,\alpha_p } = \frac{ n !_q }{ \alpha_1 !_q \cdots \alpha_p !_q } \ , 
\] 
where $ \alpha = (\alpha_1,\dots,\alpha_p) $. 
We write 
\[
 s_2(\alpha) = \sum_{i<j} \alpha_i\alpha_j 
\] 
for the second elementary symmetric function of $\alpha$. 
If $a$ is a Meshalkin sequence, we write 
\[ 
r(a) = (r(a_1),\dots,r(a_p)) 
\] 
for the sequence of ranks. 
We define $ \bbP^{\,n-1} (q) $ to be empty if $n=0$, 
a point if $n=1$, and a line of $q+1$ points if $n=2$. 

\begin{theorem}\label{main} 
Let $n\geq 0$, $l\geq 1$, $p\geq 2$, and $q\geq 2$. 
Let $\cM$ be a family of Meshalkin sequences of length $p$ in $ \bbP^{\,n-1} (q) $ 
such that, for each $k \in [p-1]$, $\cM_k$ contains no chain of length $l$. Then 
\begin{enumerate} 
\item[(a)] $ \displaystyle \sum_{a\in \cM} \dfrac{1}{\smallgauss{n}{r(a)} q^{ s_2(r(a)) }} \leq l^{p-1} $, and 
\item[(b)] $|\cM|$ is at most equal to the sum of the $l^{p-1}$ largest amongst the quantities $ \smallgauss n \alpha q^{ s_2(\alpha) } $ for $ \alpha = ( \alpha_1,\dots,\alpha_p ) $ with all $ \alpha_k \geq 0 $ and $ \alpha_1 + \dots + \alpha_p = n $. 
\end{enumerate} 
\end{theorem} 

The antichain case (where $l=1$), the analog of Meshalkin and Hochberg and Hirsch's theorems, is captured in

\begin{corollary}\label{anti} 
Let $\cM$ be a family of Meshalkin sequences of length 
$p \geq 2$ in $ \bbP^{\,n-1} (q) $ such that each $\cM_k$ for $k<p$ is an antichain. Then 
\begin{enumerate} 
\item[(a)] $ \displaystyle \sum_{a\in \cM} \dfrac{1}{\smallgauss{n}{r(a)} q^{ s_2(r(a)) }} \leq 1 $, and 
\item[(b)] $ \displaystyle |\cM| \leq \max_{\alpha} \gauss n \alpha  q^{ s_2(r(a)) } 
=  \gauss n { \big\lceil \frac np \big\rceil , \dots , \big\lceil \frac np \big\rceil , \big\lfloor \frac np \big\rfloor , \dots , \big\lfloor \frac np \big\rfloor }
   q^{s_2( \lceil n/p \rceil , \dots , \lceil n/p \rceil , \lfloor n/p \rfloor , \dots , \lfloor n/p \rfloor )} $. 
\item[(c)] Equality holds in (a) and (b) if, for each $k$, $\cM_k$ consists of all flats of rank $\big\lceil\frac n p \big\rceil$ or all of rank $\big\lfloor \frac n p \big\rfloor$.
\hfill\qedsymbol
\end{enumerate} 
\end{corollary} 

We believe---but without proof---that the largest families $\cM$ described in (c) are the only ones. 

Notice that we do not place any condition in either the theorem or its corollary on $\cM_p$. 

Our theorem is not exactly a generalization of that of Rota and Harper because a flat in a projective geometry has a variable number of complements, depending on its rank. Still, our result does imply this and a generalization, as we shall demonstrate in Section \ref{partial}. 

\section{Proof of Theorem \ref{main}}\label{proof} 

The proof of Theorem \ref{main} is adapted from the short proof of Theorem \ref{bwz} in \cite{bz}.  It is complicated by the multiplicity of complements of a flat, so we require the powerful lemma of Harper, Klain, and Rota (\cite[Lemma 3.1.3]{klain}, improving on \cite[Lemma on p.~199]{rota}; for a short proof see \cite[Lemmas 3.1 and 5.2]{bwz}) and a count of the number of complements. 

\begin{lemma}\label{hkr} 
Suppose given real numbers $ m_1 \geq m_2 \geq \dots \geq m_N \geq 0 $, 
other real numbers $ q_1,\dots,q_N \in [0,1] $, and an integer $P$ with $1 \leq P \leq N$.  If 
$ \sum_{k=1}^N q_k \leq P $, then 
\begin{equation}
\label{qm} q_1 m_1 + \dots + q_N m_N \leq m_1 + \dots + m_P \ . 
\end{equation} 

Let $m_{ P'+1 }$ and $m_{P''}$ be the first and last $m_k$'s equal to $m_P$. 
Assuming $m_P>0$, there is equality in (\ref{qm}) if and only if 
\[ 
q_k = 1 \text{ for } m_k > m_P , \qquad q_k = 0 \text{ for } m_k < m_P , \qquad \text{and} \qquad q_{P'+1} + \dots + q_{P''} = P-P' \ . 
\] 
\hfill\qedsymbol
\end{lemma} 
\smallskip

\begin{lemma}\label{comp} 
A flat of rank $k$ in $ \bbP^{\,n-1} (q) $ has $ q^{k(n-k)} $ complements. 
\end{lemma} 

\begin{proof} 
The number of ways to extend a fixed ordered basis $(P_1,\dots,P_k)$ of the flat to 
an ordered basis $(P_1,\dots,P_n)$ of $ \bbP^{\,n-1} (q) $ is 
\[ 
\frac{ q^n - q^k }{ q-1 }\ \frac{ q^n - q^{k+1} }{ q-1 }\ \cdots\ \frac{ q^n - q^{n-1} }{ q-1 } \ . 
\] 
Then $ P_{k+1} \vee \dots \vee P_n $ is a complement and is generated by the last $n-k$ points in 
\[ 
\frac{ q^{n-k} - 1 }{ q-1 }\ \frac{ q^{n-k} - q }{ q-1 }\ \cdots\ \frac{ q^{n-k} - q^{n-k-1} }{ q-1 } 
\] 
of the extended ordered bases. Dividing the former by the latter, there are 
  \[ q^{ \left( \binom n 2 - \binom k 2 \right) - \binom {n-k} 2 } = q^{k(n-k)} \]  
complements. 
\end{proof} 

\begin{proof}[Proof of (a)] 
We proceed by induction on $p$. For a flat $f$, define 
  \[ \cM(f) := \{ (a_2,\dots,a_p) : \ (f,a_2,\dots,a_p) \in \cM \} \] 
and also, letting $c$ be another flat, define 
  \[ \cM^c(f) := \{ (a_2,\dots,a_p) \in \cM(f) : \ a_2 \vee \dots \vee a_p = c \} \ . \] 
For $ a \in \cM $, we write $ r_1 = r(a_1) $. Finally, $ \cC(a_1) $ is the set 
of complements of $a_1$. If $p>2$, then 
\begin{align*} 
\sum_{a\in \cM} \dfrac{1}{\smallgauss{n}{r(a)} q^{ s_2(r(a)) }} 
                                                                           &= \sum_{a_1 \in \cM_1} \dfrac{1}{\smallgauss{n}{r_1} q^{ r_1 (n-r_1) }} \sum_{a' \in \cM (a_1) } \dfrac{1}{\smallgauss{n-r_1}{r(a')} q^{ s_2(r(a')) }} \\ 
                                                                           &= \sum_{a_1 \in \cM_1} \dfrac{1}{\smallgauss{n}{r_1} q^{ r_1 (n-r_1) }} \sum_{ c \in \cC(a_1) } \sum_{a' \in \cM^c (a_1) } \dfrac{1}{\smallgauss{n-r_1}{r(a')} q^{ s_2(r(a')) }} \\ 
                                                                           &\leq \sum_{a_1 \in \cM_1} \dfrac{1}{\smallgauss{n}{r_1} q^{ r_1 (n-r_1) }} \sum_{ c \in \cC(a_1) } l^{p-2} \\ 
\intertext{by induction, because $\cM^c(a_1)$ is a Meshalkin family in $ c \cong \bbP^{\, r(c)-1 } = \bbP^{\, n-r_1-1 } $ 
and each $\cM_k^c(a')$ for $k<p-1$, being a subset of $\cM_{k+1}$, is $l$-chain-free,} 
                                                                           &= \sum_{a_1 \in \cM_1} \dfrac{1}{\smallgauss{n}{r_1} q^{ r_1 (n-r_1) }} q^{ r_1 (n-r_1) } l^{p-2} \\ 
\intertext{by Lemma \ref{comp},} 
                                                                           &\leq l \cdot l^{p-2} \\ 
\end{align*} 
by the theorem of Rota and Harper. 

The initial case, $p=2$, is similar except that the innermost sum in the second step equals $1$. 
\end{proof}

\begin{lemma}\label{meshcount} 
Let $\alpha=(\alpha_1,\dots,\alpha_p)$ with all $\alpha_k \geq 0$ and $\alpha_1+\dots+\alpha_p=n$. 
The number of all Meshalkin sequences $a$ in $\bbP^{\,n-1}$ with $r(a)=\alpha$ is $ \smallgauss n \alpha q^{ s_2(\alpha) } $. 
\end{lemma} 

\begin{proof} 
If $p=1$, then $a = \hat 1$ so the conclusion is obvious. 
If $p>1$, we get a Meshalkin sequence of length $p$ in $\bbP^{\,n-1}$ with 
rank sequence $r(a)=\alpha$ by choosing $a_1$ to have rank $\alpha_1$, 
then a complement $c$ of $a_1$, and finally a Meshalkin sequence $a'$ of 
length $p-1$ in $ c \cong \bbP^{\, r(c)-1 } = \bbP^{\, n-\alpha_1-1 } $ 
whose rank sequence is $ \alpha' = (\alpha_2,\dots,\alpha_p)$. 
The first choice can be made in $ \gauss {n-\alpha_1}{\alpha'} $ ways, 
the second in $ q^{ \alpha_1 (n - \alpha_1) } $ ways, and the third, 
by induction, in $ \gauss{ n-\alpha_1 }{ \alpha' } q^{ s_2(\alpha') } $ 
ways.  Multiply.
\end{proof} 

\begin{proof}[Proof of (b)] 
Let $N(\alpha)$ be the number of $a \in \cM$ for which $r(a)=\alpha$.  In Lemma \ref{hkr} take 
\[ 
q_\alpha = \frac{ N(\alpha) }{ \smallgauss n \alpha q^{ s_2(\alpha) } } \qquad \text{ and } \qquad m_\alpha = \gauss n \alpha q^{ s_2(\alpha) } \ , 
\] 
and number all possible $\alpha$ so that $m_{\alpha^1} \geq m_{\alpha^2} \geq \cdots $\ .

Lemma \ref{meshcount} shows that all $ q_\alpha \leq 1 $ so Lemma \ref{hkr} does apply.  The conclusion is that 
\[ 
|\cM| = \sum_{i=1}^N q_{ \alpha^i } m_{ \alpha^i } \leq \gauss n {\alpha^1} q^{ s_2(\alpha^1) } + \dots + \gauss n {\alpha^P} q^{ s_2(\alpha^P) } \ , 
\] 
where $N = \binom{n+p-1}{p-1}$, the number of sequences $\alpha$, and $P = \min ( l^{p-1} , N )$. 
\end{proof}

\section{Strangeness of the LYM Inequality}\label{strange} 

There is something odd about the LYM inequality in Theorem \ref{main}(a). 
A normal LYM inequality would be expected to have denominator $ \smallgauss n {r(a)} $ 
without the extra factor $q^{ s_2(r(a)) } $. Such an LYM inequality does exist; it is a 
corollary of Theorem \ref{main}(a); but it is not strong enough to give the upper 
bound on $|\cM|$. We prove this weaker inequality here. 

\begin{proposition}\label{weaklym} 
Assume the hypotheses of Theorem \ref{main}; that is: 
$n\geq 0$, $l\geq 1$, $p\geq 2$, and $q\geq 2$; and 
$\cM$ is a family of Meshalkin sequences of length $p$ in $ \bbP^{\,n-1} (q) $ 
such that, for each $k \in [p-1]$, $\cM_k$ contains no chain of length $l$. Then 
$ \sum_{ a \in \cM } \smallgauss n {r(a)}^{ -1 } $ is bounded above by the sum of 
the $l^{p-1}$ largest expressions $q^{s_2(\alpha)}$ for $ \alpha = ( \alpha_1 , \dots , \alpha_p ) $ 
with all $\alpha_k \geq 0$ and $\alpha_1 + \dots + \alpha_p = n$. 
\end{proposition} 
\begin{proof} 
Again we apply Lemma \ref{hkr}, this time with $ q_\alpha = {N(\alpha)}/{ \gauss n \alpha q^{ s_2(\alpha) } } $ 
and $ M_\alpha = q^{ s_2(\alpha) } $. 
\end{proof} 

\section{A ``Partial" Corollary}\label{partial} 

We deduce Theorem \ref{rh}(a) from the case $p=2$ of Theorem \ref{main}(a). 
Our purpose is not to give a new proof of Theorem \ref{rh} but to show that
we have a generalization of it. 

The key to the proof is that $\cM_2$ in our theorem is not required to be 
$l$-chain-free. Therefore if we have an $l$-chain-free set $\cA$ of flats in 
$\bbP^{\,n-1}$, we can define 
\[ 
\cM = \left\{ (a,c) : \ a \in \cA \text{ and } c \in \cC(a) \right\} \ ; 
\] 
and $\cM$ will satisfy the requirements of Theorem \ref{main}. 
The LYM sum in Theorem \ref{main}(a) then equals the LYM sum 
in Theorem \ref{rh}(a), and we are done. 

The same argument gives a general corollary. A \emph{partial Meshalkin sequence of length $p$} 
is a sequence $ a = (a_1,\dots,a_p) $ of flats in $\bbP^{\,n-1}(q)$ such that 
$ r(a_1 \vee \dots \vee a_p) = r(a_1)+ \dots + r(a_p) $. We simply do not require 
the join $ \hat a = a_1 \vee \dots \vee a_p $ to be $ \hat 1 $. The generalized 
Rota--Harper theorem is: 

\begin{corollary}\label{partialmain} 
Let $p \geq 1$, $l \geq 1$, $q \geq 2$, and $n \geq 0$. Let $\cM$ be a family of 
partial Meshalkin sequences of length $p$ 
in $\bbP^{\,n-1}(q)$ such that, for each $k \in [p]$, $\cM_k$ contains no chain of length $l$. Then 
\begin{enumerate} 
\item[(a)] $ \displaystyle \sum_{ a \in \cM } \frac 1 { \smallgauss n { r ( \hat a ) } \smallgauss {r ( \hat a )}{r(a)} q^{ s_2 (r(a)) } } \leq l^p $ and 
\item[(b)] $|\cM|$ is at most equal to the sum of the $l^p$ largest amongst the quantities $ \gauss n \alpha q^{ s_2(\alpha) } $ for $ \alpha = ( \alpha_1 , \dots , \alpha_{p+1} ) $ with all $\alpha_k \geq 0 $ and $ \alpha_1 + \dots + \alpha_{p+1} = n $.
\hfill\qedsymbol 
\end{enumerate} 
\end{corollary} 

As a special case we generalize the $q$-analog of Sperner's theorem. (The $q$-analog is the case $p=1$.) 

\begin{corollary}\label{partialanti} 
Let $\cM$ be a family of partial Meshalkin sequences of length $p \geq 1$ in $\bbP^{\,n-1}$ 
such that each $\cM_k$ is an antichain. Then: 
\begin{enumerate} 
\item[(a)] $ \displaystyle \sum_{ a \in \cM } \frac 1 { \smallgauss n { r ( \hat a ) } \smallgauss {r ( \hat a )}{r(a)} q^{ s_2 (r(a)) } } \leq 1 $. 
\item[(b)] $|\cM| \leq \gauss {n}{\alpha}q^{s_2(\alpha)}$, in which $\alpha = \big( \big\lceil \frac n {p+1} \big\rceil , \dots , \big\lceil \frac n {p+1} \big\rceil , \big\lfloor \frac n {p+1} \big\rfloor , \dots , \big\lfloor \frac n {p+1} \big\rfloor  \big)$ where the number of terms equal to $ \big\lceil \frac n {p+1} \big\rceil $ is the least nonnegative residue of $n$ modulo $p+1$. 
\item[(c)] Equality holds in (a) and (b) if, for each $k$, $\cM_k$ consists of all flats of rank $ \big\lceil \frac n {p+1} \big\rceil $ or all flats of rank $ \big\lfloor \frac n {p+1} \big\rfloor $. 
\hfill\qedsymbol
\end{enumerate} 
\end{corollary} 

We conjecture that the largest families $\cM$ described in (c) are unique. 

% {\bf Acknowledgements}. 

% -------------------------------------------

%\newpage
%\small
%\footnotesize
%\nocite{*}
%\addcontentsline{toc}{subsubsection}{References}
\bibliography{thesis}
\bibliographystyle{alpha}

\end{document}